\documentclass[A4, 10 pt, conference]{ieeeconf}  

\IEEEoverridecommandlockouts                              

\usepackage{amsmath} 
\usepackage{amssymb}  
\usepackage{bm}
\usepackage{amstext}
\usepackage{mathtools}
\usepackage{tikz}
\usepackage{cite}
\usepackage{enumerate}
\usetikzlibrary{shapes,arrows,circuits,decorations.pathmorphing,patterns,positioning,decorations.markings,calc}
\usepackage{pgfplots}
\tikzstyle{block} = [draw, fill=blue!20, rectangle, minimum height=3em, minimum width=6em]
\tikzstyle{sum} = [draw, fill=blue!20, circle, node distance=1cm]
\tikzstyle{input} = [coordinate]
\tikzstyle{output} = [coordinate]
\tikzstyle{pinstyle} = [pin edge={to-,thin,black}]

\newtheorem{definition}{Definition}[section]
\newtheorem{theorem}{Theorem}[section]

\newtheorem{lemma}[theorem]{Lemma}

\newtheorem{assumption}{Assumption}
\newtheorem{remark}{Remark}

\title{\LARGE \bf
On the achievable degree of stability in strictly negative imaginary state feedback control
}

\author{James Dannatt$^{1}$, Ian Petersen$^{2}$
	\thanks{*This work was supported by the Australian Research Council under grants DP160101121 and DP190102158.}
	\thanks{$^{1}$James Dannatt is with the Research School of Electrical, Energy and Materials Engineering at the Australian National University, ACT.
		{\tt\small james.dannatt@anu.edu.au}.}%
	\thanks{$^{2}$Ian Petersen is with the Research School of Electrical, Energy and Materials Engineering at the Australian National University, ACT.
		{\tt\small i.r.petersen@gmail.com}.}%
}

\begin{document}

\maketitle
\thispagestyle{empty}
\pagestyle{empty}

\begin{abstract}
In this paper we study the problem of determining the largest degree of stability that can be achieved for SISO systems using negative imaginary state feedback control. A state feedback result is given for synthesising a controller for a plant such that a given closed-loop transfer function is strictly negative imaginary with a prescribed degree of stability. By varying a parameter, the placement of the closed-loop system poles can be adjusted to give a prescribed degree of stability. We show the achievable degree of stability is related to the zero locations of the transfer function from the control input to the disturbance output of the nominal plant being controlled. Moreover, we offer results that outline the largest degree of stability that can be achieved for systems with distinct eigenvalues.
\end{abstract}

\section{Introduction}
The theory of negative imaginary (NI) systems is broadly applicable to problems of robust vibration control for flexible structures; e.g., see \cite{LP06,LP_CSM,CH10a,MKPL11}. In these control systems, unmodelled spillover dynamics can degrade control system performance or lead to instability if the controller is not designed to be robust against this type of uncertainty. NI systems theory provides a way analyzing robustness and designing robust controllers for such flexible structures in the case of collocated force actuators and position sensors; e.g., see \cite{BMP1,CH10a,MKPL1,MLKP1,MKPL5,MKPL11,XiLaP1}. Motivated by the robust stability properties of NI systems, research has been made into controller synthesis results with the aim of creating a closed-loop system with the NI or SNI property \cite{LP06,XiPL1}.

\vspace{10pt}

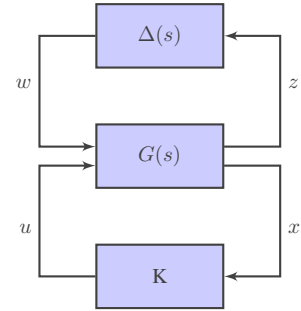
\begin{figure}[!ht]
	\centering
	\begin{tikzpicture}[auto, node distance=2cm,>=latex', black!75,thick, scale=0.8, every node/.style={transform shape}]
		\node [block, name=uncertainty] {$\Delta(s)$};
		\node [block, below of=uncertainty] (plant) {$G(s)$};
		\node [block, below of=plant] (controller) {K};
		\node [coordinate, left of=plant]  (input) {};
		\node [coordinate, right of=plant] (output) {};
		
		\draw [->] (uncertainty) -| node[left, yshift=-5ex]{$w$}([yshift= 1ex]input) -- ([yshift= 1ex]plant.west) ;
		\draw [<-] ([yshift= -1ex]plant.west) -- ([yshift= -1ex]input) |- node[left, yshift=5ex]{$u$} (controller);
		\draw [->] ([yshift= 1ex]plant.east) -- ([yshift= 1ex]output) |- node[right, yshift=-5ex]{$z$} (uncertainty);
		\draw [<-] (controller) -| node[right, yshift=5ex]{$x$}([yshift= -1ex]output) -- ([yshift= -1ex]plant.east);
	\end{tikzpicture}
	\caption{A state feedback control system with plant uncertainty $\Delta(s)$.}
	\label{F1}
\end{figure}

In this paper, we concentrate on the state feedback negative imaginary control problem as illustrated in Figure \ref{F1}. In this control problem, the unmodelled flexible dynamics are represented by the plant uncertainty $\Delta(s)$ which are assumed to have the NI or SNI property. This is a natural assumption in the control of flexible structures, since it is known that any flexible structure with collocated force inputs and position outputs will have the NI property; e.g., see \cite{LP_CSM}. The (strictly) NI state feedback control problem is concerned with synthesising a state feedback control law $u = Kx$, such that the corresponding closed-loop system will be guaranteed to have the (strictly) NI property

One method of approaching the SNI state feedback control problem has been to form a controller using the solution to a certain algebraic Riccati equation (ARE) \cite{salcan2018negative,auto_reyes_2019}. However, it has been shown that the associated ARE used in these results does not have a stabilizing solution \cite{DaPL1}. To avoid the computation complexity of solving this ARE, \cite{MKPL6a} proposed a controller synthesis approach using the solution to an algebraic Riccati equation (ARE) that could be obtained by solving two Lyapunov equations. The approach was computationally efficient, but yielded a closed-loop pole at the origin, ensuring a marginally stable closed-loop system. This pole at the origin was problematic, as in vibration control problems for flexible structures, pole placement is related to the degree of damping achieved for the nominal resonant modes. 

\vspace{10pt}

The papers \cite{MKPL7a,MKPL11} modified the approach of \cite{MKPL6a} using a perturbation applied to the plant matrix of the nominal system in order to ensure asymptotic stability of the closed-loop system. The perturbation achieved closed-loop asymptotic stability. However, it was not proven that the closed-loop system retained the NI property after such a perturbation. The paper \cite{DaPL1} then showed that under suitable assumptions, the perturbation approach of \cite{MKPL11} did in fact lead to an SNI closed-loop system and by varying the perturbation parameter, the closed-loop poles of the system could be shifted further into the complex plane, resulting in a prescribed degree of stability. The Lyapunov equations used in forming the controller in this approach also depend on the perturbation parameter. As a consequence, there are choices of this parameter that result in the equations no longer being solvable. This puts a limit on the degree of stability that can be achieved.

\vspace{10pt}

In this paper, we focus on SISO systems and build on the results of \cite{DaPL1} to study the achievable degree of stability of an SNI closed-loop system formed using the perturbation approach of \cite{MKPL11}. This issue is important since closed-loop degree of stability is a key performance measure in vibration control problems for flexible structures. Our main results show that the achievable degree of stability for a given system is limited by the eigenvalues of a certain matrix that is related to the zeros of the transfer function from the control input to the disturbance output of the nominal system. Under suitable assumptions, we quantify how the perturbation approach of \cite{MKPL11} effects this matrix, and the maximum degree of stability that can be achieved before we can no longer guarantee the SNI property of the closed-loop system.

\section{Preliminary Definitions}
In this section, we briefly present definitions for both negative imaginary and strictly negative imaginary systems.

First consider the linear time-invariant (LTI) system,
\begin{eqnarray} \label{system}
	\dot{x}(t) &=& Ax(t) + Bu(t), \nonumber \\
	y(t) &=& Cx(t) + Du(t),
\end{eqnarray}
where $A \in \mathbb{R}^{n \times n},$ $B \in \mathbb{R}^{n \times m},$ $C \in \mathbb{R}^{m \times n}$ and $D \in \mathbb{R}^{m \times m}$.

The following two definitions relate to the NI and SNI properties of the transfer function matrix $G(s) = C(sI - A)^{-1}B+D$ corresponding to the system (\ref{system}). We do not consider a direct feed-through term from disturbance to output in this work.

\vspace{10pt}
\begin{definition} \label{def:NI}{ A square transfer function matrix $G(s)$ is NI if the following conditions are satisfied \cite{MKPL10}:}
	\begin{enumerate}[(i)]
		\item G(s) has no pole in $Re[s]>0$.

		\item For all $\omega \geq 0$ such that $jw$ is not a pole of $G(s)$, $j(G(j\omega) - G(j\omega)^*) \geq 0$.

		\item If $s=j\omega_0$, $ \omega_0 > 0$ is a pole of $G(s)$ then it is a simple pole. Furthermore, if $s=j\omega_0$, $ \omega_0 > 0$ is a pole of $G(s)$, then the residual matrix $K = \lim_{s \to j\omega_0} (s-j\omega_0)jG(s)$ is positive semidefinite Hermitian.

		\item If $s=0$ is a pole of $G(s)$, then it is either a simple pole or a double pole. If it is a double pole, then, $\lim_{s \to 0} s^2G(s) \geq 0$.
	\end{enumerate}
\end{definition}
\vspace{10pt}
Also, an LTI system (\ref{system}) is said to be NI if the corresponding transfer function matrix $G(s) = C(sI-A)^{-1}B + D$ is NI.
\vspace{10pt}
\begin{definition} \label{def: SNI freq domain definition}
{A square transfer function matrix $G(s)$ is SNI if the following conditions are satisfied \cite{MKPL10}:}
	\begin{enumerate}[(i)]
		\item G(s) has no poles in $Re[s] \geq 0$.

		\item For all $\omega > 0$ such that $jw$ is not a pole of $G(s)$, $j(G(j\omega) - G(j\omega)^*) > 0$.
	\end{enumerate}
\end{definition}

Also, an LTI system (\ref{system}) is said to be SNI if the corresponding transfer function matrix $G(s) = C(sI-A)^{-1}B + D$ is SNI.
\vspace{10pt}
\begin{lemma}(\cite{REY19}) \label{ch: NIrealisation: lemma: non-minimal NI ARE lemma} Suppose the system (\ref{system}) is a given state space realization with $R = CB + B^TC^T > 0$ and $D=D^T$. If there exists a real $P=P^T\geq0$ such that
	\begin{equation} \label{chap NIrel: NI ARE WITH semi-stable A}
		PA + A^TP + (CA - B^TP)^TR^{-1}(CA-B^TP) = 0
	\end{equation}
	with $\sigma(A-BR^{-1}(CA-B^TP))\subset \mathbb{C}_{\leq 0}$. Then, the state space realization and corresponding transfer function is NI.
\end{lemma}

\section{Strictly Negative Imaginary State Feedback}

In this paper are concerned with the use of state feedback control to achieve a closed-loop system which is SNI and has a prescribed degree of stability. Moreover, we are concerned with how stable can we make the eigenvalues of the closed-loop system while maintaining the SNI property. In this section, we introduce the SNI state feedback control problem along with the state feedback results of the papers \cite{MKPL11,DaPL1}. In addition to this we  will introduce supporting lemmas that contextualise the perturbation method used in the state feedback theorem of \cite{MKPL11,DaPL1}.

\subsection{Problem Formulation}

For the SNI state feedback method we will consider in this paper, we require the following controlled linear state space system:
\begin{align}
	\dot{x} &= Ax + B_1w + B_2u, \nonumber\\
	z       &= C_1x, \label{uncertain system}
\end{align}
where $A \in \mathbb{R}^{n \times n}$, $B_1 \in \mathbb{R}^{n \times m}$, $B_2 \in \mathbb{R}^{n \times m}$, and $C_1 \in \mathbb{R}^{m \times n}$. This system corresponds to the nominal plant transfer function matrix $G(s)$ in Figure \ref{F1}.

If we apply a state feedback control law $u=Kx$ to this system, the corresponding closed-loop system has state space representation
\begin{align}
	\dot{x} &= (A + B_2K)x + B_1w, \nonumber\\
	z       &= C_1x,        \label{math: closed-loop system}
\end{align}
with corresponding closed-loop transfer function $G_{cl}(s) = C_1(sI - A - B_2K)^{-1}B_1$. 

\vspace{10pt}

We assume that the system (\ref{uncertain system}) satisfies the following assumptions.
\begin{assumption}
  \label{A1}
  The matrix $C_1B_2$ is non-singular.
\end{assumption}
\begin{assumption}
  \label{A2}
  The matrix $R = C_1B_1 + B_1^TC_1^T > 0$.	
\end{assumption}
Assumption \ref{A1} is essential for ensuring the existence of the controller in the state feedback theorem that will follow. Assumption \ref{A2} is necessary to ensure the invertability of a certain matrix in several of the proofs that follow.

\subsection{Preliminary Results}

We now introduce three results that are needed to understand the SNI state feedback theorem discussed in this paper. Firs,t we will introduce a lemma that details the behavior of NI systems under an $\epsilon{I}$ perturbation applied to the plant matrix. Then, we present a Schur decomposition applied to a matrix that is used in controller synthesis to avoid the potential computational difficulties of solving an ARE that does not have a stabilizing solution. In the third result, we remark that the eigenvalues of the matrix used in the Schur decomposition actually corresponds to the zeros of a certain related transfer function. Our main results will show that these eigenvalues dictate the maximum degree of stability that the closed-loop system can achieve under SNI state feedback.

\vspace{10pt}

\begin{lemma}[\cite{DaPL1}] \label{lemma: NI purturb results in SNI} Suppose the system (\ref{system}) is a given state space realization with $m=1$ and an NI transfer function $G(s)$. Then the perturbed state space realization {\tiny$\begin{bmatrix}
			\begin{tabular}{ l | r }
				$A-\epsilon I$ & $B$ \\ \hline
				$C$ & $D$
			\end{tabular}
		\end{bmatrix}$} will be SNI for all $\epsilon > 0$.
\end{lemma}

\vspace{10pt}

We now introduce the matrix
\begin{align} 
	A_r &= Q(A+\epsilon I) = A_q + \epsilon Q.  \label{A_r equation} 
\end{align}
where 
\begin{align}
  Q &= I-B_2(C_1B_2)^{-1}C_1, \label{math: schur Q definition}\\
  	A_q &= QA, \label{math: A_q equation}
\end{align}
and $\epsilon > 0$ is a parameter which will determine the degree of stability of the closed-loop system in the SNI state feedback theorem of the next subsection. The matrix $A_r$ is fundamental to the state feedback theorem that follows and is discussed further in Remark~\ref{R1}.

\vspace{10pt}

After applying a real Schur decomposition to the matrix $A_r$ as in (\ref{chap how far pert:math: AF}), we may take the orthogonal matrix $U$ and construct the following matrices:
\begin{subequations} \label{math: Schur decomp equations}
	\begin{align}
		\tilde{A} &= U^TA_rU =
		\begin{bmatrix}
			\begin{tabular}{ l  r }
				$\tilde{A}_{11}$ & $\tilde{A}_{12}$ \\
				0 & $\tilde{A}_{22}$
			\end{tabular}
		\end{bmatrix}, \label{chap how far pert:math: AF}\\
		\tilde{B} &= U^TB_1 =
		\begin{bmatrix}
			\begin{tabular}{c}
				$\tilde{B}_{11}$ \\
				$\tilde{B}_{22}$
			\end{tabular}
		\end{bmatrix}, \label{chap how far pert:math: B1}\\
		\tilde{C} &= U^T\big(B_2(C_1B_2)^{-1} - B_1R^{-1}\big) =
		\begin{bmatrix}
			\begin{tabular}{c}
				$\tilde{C}_{11}$ \\
				$\tilde{C}_{22}$
			\end{tabular}
		\end{bmatrix}, \label{chap how far pert:math: Bf}\\
		\tilde{Z} &= U^TZU = \tilde{B}R^{-1}\tilde{B}^T - \tilde{C}R\tilde{C}^T = \begin{bmatrix}
			\begin{tabular}{ l  r }
				$\tilde{Z}_{11}$ & $\tilde{Z}_{12}$ \\
				$\tilde{Z}_{21}$ & $\tilde{Z}_{22}$
			\end{tabular}
		\end{bmatrix}, \label{chap how far pert:math: Zf}
	\end{align}
\end{subequations}
where $\tilde{A}_{11}$ has all of its eigenvalues in the closed left half of the complex plane and $\tilde{A}_{22}$ is an anti-stable matrix; i.e., $\sigma(\tilde{A}_{22}) \subset \{s:\operatorname{Re}[s]>0\}$. Here, $U$ is an orthogonal matrix dependent on $\epsilon$ and is obtained through the real Schur transformation; see Section 5.4 of \cite{BER09}. Also, 
\begin{align}
	Z &= B_1(B_2^TC_1^T)^{-1}B_2^T + B_2(C_1B_2)^{-1}B_1^T \nonumber \\
	& \qquad \qquad \qquad - B_2(C_1B_2)^{-1}R(B_2^TC_1^T)^{-1}B_2^T. \label{math: schur Z definition}
\end{align} 

\vspace{10pt}

\begin{remark}
  \label{R1}
Note that the choice of matrix $A_r$ for the Schur decomposition above is not arbitrary. The eigenvalues of $A_q$ correspond to the zeros of the transfer function from $u$ to $z$ in (\ref{uncertain system}), plus a zero at the origin. To see this, set $w\equiv 0$ and note that (\ref{uncertain system}) implies
\begin{align}
	\dot{z} = C_1Ax+C_1B_2u. \label{pertb ch: math z dot}
\end{align}
Since $(C_1B_2)$ is invertible by assumption, we can rearrange (\ref{pertb ch: math z dot}) to write
\begin{align}
	u = (C_1B_2)^{-1}\dot{z} - (C_1B_2)^{-1}C_1Ax. \label{pertb ch: math u as function of zdot}
\end{align}
Now, we can substitute (\ref{pertb ch: math u as function of zdot}) into (\ref{uncertain system}) to give
\begin{align}
	\dot{x} &= Ax + B_2(C_1B_2)^{-1}\dot{z}- B_2(C_1B_2)^{-1}C_1Ax \nonumber\\
	&= (A- B_2(C_1B_2)^{-1}C_1A)x + B_2(C_1B_2)^{-1}\dot{z} \nonumber\\
	&= A_qx + B_2(C_1B_2)^{-1}\dot{z}. \label{pertb ch: math x dot}
\end{align}
These equations define the inverse system which maps from $\dot{z}$ to $u$. Considering the Laplace transform of $\dot{z}$, suppose $z(0)=0$ and observe that $\mathcal{L}(\dot{z}) = s\mathcal{L}(z) + z(0)$. Thus, the eigenvalues of $A_q$ will be the zeros of the transfer function from $u$ to $z$ plus a zero at zero.
\end{remark}

\subsection{Negative Imaginary State Feedback}

In this section we introduce the SNI state feedback theorem that we will be analysing in our main result to determine the degree of stability that can be achieved. This theorem can be found in \cite{MKPL11,DaPL1}.

\vspace{10pt}

\begin{theorem}[\cite{DaPL1}]\label{theorem: SNI ARE synthesis theorem} Consider the  system (\ref{uncertain system}) with $m=1$ satisfying Assumptions A1-A2. For a given $\epsilon > 0$, there exists a static state feedback matrix $K$ such that the closed-loop system (\ref{math: closed-loop system}) is SNI with degree of stability $\epsilon$ if there exist  $T \geq 0$ and $S \geq 0$ such that
	\begin{align}
		-\tilde{A}_{22}T - T\tilde{A}_{22}^T + \tilde{C}_{22}R\tilde{C}_{22}^T &= 0, \label{math: T cond}\\
		-\tilde{A}_{22}S - S\tilde{A}_{22}^T + \tilde{B}_{22}R^{-1}\tilde{B}_{22}^T &= 0, \label{math: S cond}\\
		T-S &> 0, \label{math: T-S cond}
	\end{align}
	where $\tilde{A}_{22}$, $\tilde{B}_{22}$ and $\tilde{C}_{22}$ are obtained from the Schur decomposition (\ref{math: Schur decomp equations}). Moreover, if the conditions (\ref{math: T cond})-(\ref{math: T-S cond}) are satisfied, then the required state feedback controller matrix $K$ is given by
	\begin{equation} 
		K = (C_1B_2)^{-1}(B_1^TP - C_1A - \epsilon C_1 - R(B_2^TC_1^T)^{-1}B_2^TP), \label{math: k matrix}
	\end{equation} where $P = U\tilde{P}U^T$ and $\tilde{P} = \begin{bmatrix}
		\begin{tabular}{ l  r }
			$0$ & $0$ \\
			$0$ & $(T-S)^{-1}$
		\end{tabular}
	\end{bmatrix} \geq 0$. Here, $U$ is the orthogonal matrix obtained through the Schur transformation (\ref{math: Schur decomp equations}).
\end{theorem}

\vspace{10pt}

When the $A_q$ matrix from the Schur decomposition (\ref{math: Schur decomp equations}) has no unstable eigenvalues, the dimension of $\tilde{A}_{22}$ is equal to zero. Therefore, we cannot form (\ref{math: T cond})-(\ref{math: T-S cond}). This situation corresponds to the case in which the transfer function from $u$ to $y$ in (\ref{uncertain system}) is minimum phase. This case is also considered in the paper \cite{SPV4a} for the special case of $B_1= B_2$ in (\ref{uncertain system}). The following theorem is a minor result that allows us to use SNI state feedback in the case in which the transfer function from $u$ to $y$ in (\ref{uncertain system}) is minimum phase.

\vspace{10pt}

\begin{theorem}\label{theorem: SNI ARE synthesis theorem no A22} Consider the  system (\ref{uncertain system}) with $m=1$ that satisfies Assumptions A1-A2. Suppose the matrix $A_q$, obtained from the Schur decomposition (\ref{math: Schur decomp equations}) has no unstable eigenvalues. For a given $\epsilon > 0$ such that the corresponding matrix $A_r$ also has no unstable eigenvalues, there exists a static state feedback matrix $K$ such that the closed-loop system (\ref{math: closed-loop system}) is SNI with degree of stability $\epsilon$. The required state feedback controller matrix $K$ is given by
	\begin{equation}
		K = -(C_1B_2)^{-1}(C_1A + \epsilon C_1 ). \label{math: k matrix reduced}
	\end{equation} 
\end{theorem}

\vspace{10pt}

{\em Proof of Theorem \ref{theorem: SNI ARE synthesis theorem no A22}:} 
	For design purposes we will begin by replacing the open loop plant matrix $A$ with the perturbed matrix $A_\epsilon = A + \epsilon{I}$. Now, using Lemma~\ref{ch: NIrealisation: lemma: non-minimal NI ARE lemma}, we know that the closed-loop system (\ref{math: closed-loop system}) is NI if there exists a real $P=P^T\geq0$ such that
	\begin{equation}
		P\hat{A} + \hat{A}^TP + PB_1R^{-1}B_1^TP + Q = 0 \label{math: state feedback no a22 ARE}
	\end{equation}
	where
	\begin{align*}
		\hat{A} &= A_{cl} - B_1R^{-1}C_1A_{cl}, &
		R & = C_1B_1 + B_1^TC_1^T, \\
		Q &= A_{cl}^TC_1^TR^{-1}C_1A_{cl}, &
		A_{cl} &= A_\epsilon + B_2K,
	\end{align*}
	with $\sigma \bigg(A_{cl}-B_1R^{-1}(CA_{cl}-B^TP)\bigg) \subset \mathbb{C}_{\leq 0}$. 
	
	Let $P=0$. Then, (\ref{math: state feedback no a22 ARE}) reduces to the condition
	\begin{align}
		A_{cl}^TC_1^TR^{-1}C_1A_{cl} = 0.
	\end{align}
	We will consider the matrix $C_1A_{cl}$. For our choice of $K$, and noting that for a SISO system, $(C_1B_2)$ is a scalar, $C_1A_{cl} = A+\epsilon I + B_2K$ reduces to
	\begin{align}
		C_1A_{cl} &= C_1A+\epsilon C_1 -C_1B_2(C_1B_2)^{-1}(C_1A + \epsilon C_1) = 0 \nonumber\\
	\end{align}
	Therefore, $P=0$ satisfies (\ref{math: state feedback no a22 ARE}). We will now show that $\sigma \bigg(A_{cl}-B_1R^{-1}(CA_{cl}-B^TP)\bigg) \subset \mathbb{C}_{\leq 0}$. We begin by considering $A_{cl}-B_1R^{-1}(CA_{cl}-B_1^TP)$. When $P=0$, and substituting for $K = -(C_1B_2)^{-1}(C_1A + \epsilon C_1 )$ in this expression, we see that
	\begin{align}
		A_{cl}-B_1R^{-1}(CA_{cl}-B_1^TP) =& A_{cl}-B_1R^{-1}CA_{cl} \nonumber\\
		=& A_\epsilon - B_2(C_1B_2)^{-1}C_1A_\epsilon \nonumber\\
		=& A_r, \label{math: acl = ar}
	\end{align}
	where $A_r$ is defined in (\ref{A_r equation}).
	Thus, (\ref{math: acl = ar}) implies
	\begin{align}
		\sigma(A_{cl}-BR^{-1}(CA_{cl}-B^TP)) =& \sigma(A_r).
	\end{align}
	Since $\sigma(A_r)\subset \mathbb{C}_{\leq 0}$ follows from our assumption that $A_r$ has no unstable eigenvalues, we have $\sigma(A_{cl}-BR^{-1}(CA_{cl}-B^TP)) \subset \mathbb{C}_{\leq 0}$. Therefore, it follows from Lemma~\ref{ch: NIrealisation: lemma: non-minimal NI ARE lemma} that the closed-loop system is NI. Further, Lemma~\ref{lemma: NI purturb results in SNI} then implies that the actual closed-loop system corresponding to the unperturbed system will have all its poles shifted by an amount $\epsilon$ to the left in the complex plane. Moreover, the actual closed-loop system is SNI with degree of stability $\epsilon$.
\hfill $\blacksquare$

\vspace{10pt}

\section{The Achievable Degree of Stability}

This section contains the main results of this paper and provides a solution to the maximum degree of stability that can be achieved using Theorem~\ref{theorem: SNI ARE synthesis theorem} to design a controller for SNI state feedback. This section is broken into two subsections. The first offers preliminary results needed in the proofs of the main results. In the second subsection, we offer two results that outline the maximum degree of stability that can be achieved using Theorem~\ref{theorem: SNI ARE synthesis theorem}. There are two results because depending on the zeros of the transfer function from $u$ to $z$ in (\ref{uncertain system}), the matrix $A_r$ may have no unstable poles. As a consequence, Theorem~\ref{theorem: SNI ARE synthesis theorem} cannot be used and Theorem~\ref{theorem: SNI ARE synthesis theorem no A22} must be evaluated instead.

\subsection{Preliminary Results}

In this section we provide lemmas necessary to the proof of our main results.

\begin{lemma} \label{chap how far pert:Lemma: C1 orthogonal to w}
	Consider the  system (\ref{uncertain system}) and suppose the dimensions are such that $(C_1B_2) \neq 0$ is a scalar. In addition, let
	\begin{align} \label{equation: w} 
		{w} = \frac{(C_1B_1B_2^T-C_1B_2B_1^T)}{(C_1B_2)^2}.
	\end{align}
	The vector $C_1$ is orthogonal to ${w}$.
\end{lemma}

\begin{proof}
	We show our result using a simple algebraic manipulation. Observe that
	\begin{align}
		{w} C_1^T &= (\frac{C_1B_1B_2^T-C_1B_2B_1^T}{(C_1B_2)^2})C_1^T \nonumber\\
		&= \frac{C_1B_1B_2^TC_1^T-C_1B_2B_1^TC_1^T}{(C_1B_2)^2} \nonumber\\
		&= 0. \label{chap how far pert:math: w^tC_1^T = 0 equation}
	\end{align}
\end{proof}

\vspace{10pt}

\begin{lemma} \label{chap how far pert: statefeedback chap: SNI sythn: C_1 common to Aq Ar and Q}
	The vector $C_1$ defined in (\ref{uncertain system}) is in the left-nullspace of the matrices $Q$, $A_q$ and $A_r$, where $Q$, $A_q$ and $A_r$ are defined in (\ref{math: schur Q definition}), (\ref{math: A_q equation}) and (\ref{A_r equation}).
\end{lemma}
\begin{proof} 
	We can show that $C_1$ is in the left-nullspace of $Q$, $A_q$ and $A_r$ using the following routine algebra:
	\begin{align}
		C_1Q &= C_1(I-B_2(C_1B_2)^{-1}C_1) \nonumber\\
		&= C_1-C_1B_2(C_1B_2)^{-1}C_1 = C_1-C_1 = 0, \\
		C_1A_q &= C_1QA = 0, \\	
		C_1A_r &= C_1(QA + \epsilon{Q}) = C_1QA + C_1Q\epsilon = 0 + 0 = 0.
	\end{align} 
\end{proof}

\vspace{10pt}

\begin{lemma} \label{chap how far pert:Lemma: Z cases for X pos def} 
	Consider the  system (\ref{uncertain system}) with $m=1$ satisfying Assumptions A1-A2. Let $Z$ be defined as in (\ref{math: schur Z definition}). For some non-zero vector ${y}$, ${y}^T Z {y}=0$ if and only if
	\begin{align} \label{chap how far pert:conditions: SNI two conditions for X > 0}
		\frac{(C_1B_1B_2^T-C_1B_2B_1^T)}{(C_1B_2)^2}{y} &= 0, \\
		\text{or  } \qquad  B_2^T{y} &= 0. \label{chap how far pert:conditions: SNI two conditions for X > 0, cond b}
	\end{align}
\end{lemma}

\begin{proof}
	Note that Assumption A1 is required to ensure $C_1B_2 \neq 0$ in the manipulation to follow. Now, suppose there exists a non-zero vector ${y}$ which we use to pre and post multiply $Z$ to form ${y}^TZ{y}$. We may expand ${y}^TZ{y}=0$ as
	\begin{align}
		0 &= {y}^T\bigg[B_2(C_1B_2)^{-1}R(B_2^TC_1^T)^{-1}B_2^T - B_1(B_2^TC_1^T)^{-1}B_2^T \nonumber\\
		&- B_2(C_1B_2)^{-1}B_1^T\bigg]{y}. \label{chap how far pert:math: z partial equation}\\
		\intertext{Using that fact that $R = C_1B_1 + B_1^TC_1^T = 2C_1B_1$ and $C_1B_2$ are scalars, we may further manipulate (\ref{chap how far pert:math: z partial equation}) as follows}
		0 &= \frac{{y}^TB_2(2C_1B_1)B_2^T{y}}{(C_1B_2)^2} - \frac{{y}^TB_1B_2^T{y}}{C_1B_2} - \frac{{y}^TB_2B_1^T{y}}{C_1B_2} \nonumber\\
		&= \frac{(2C_1B_1)B_2^T{y}B_2^T{y}}{(C_1B_2)^2} - \frac{B_1^T{y}B_2^T{y}}{C_1B_2} - \frac{B_1^T{y}B_2^T{y}}{C_1B_2} \nonumber\\
		&= (\frac{(2C_1B_1B_2^T){y}}{(C_1B_2)^2} - \frac{B_1^T{y}}{C_1B_2} - \frac{B_1^T{y}}{C_1B_2})B_2^T{y} \nonumber\\
		&= 2(\frac{C_1B_1B_2^T}{(C_1B_2)^2} - \frac{B_1^T}{(C_1B_2)}){y}B_2^T{y} \nonumber\\
		&= 2\frac{(C_1B_1B_2^T-C_1B_2B_1^T)}{(C_1B_2)^2}{y}B_2^T{y}. \label{chap how far pert:math: when is Z zero?} 
	\end{align}
	Now, there are only two cases where (\ref{chap how far pert:math: when is Z zero?}) is equal to zero, either:
	\begin{align} 
		\frac{(C_1B_1B_2^T-C_1B_2B_1^T)}{(C_1B_2)^2}{y} &= 0, \\
		\text{or  } \qquad   B_2^T{y} &= 0,
	\end{align}
	as required.
\end{proof}

\vspace{10pt}

\begin{lemma} \label{ch how far perturb: lem Z<0 implies X > 0 : non scalar} 
	Consider the  system (\ref{uncertain system}) with $m=1$ satisfying Assumptions A1-A2. Let $\epsilon > 0$ be given and $A_r$ defined as in (\ref{A_r equation}) and suppose there exist matrices $T, S \geq 0$ satisfying:
	\begin{align}
		-\tilde{A}_{22}T - T\tilde{A}_{22}^T + \tilde{C}_{22}R\tilde{C}_{22}^T &= 0, \label{chap how far pert:math: Ts cond}\\
		-\tilde{A}_{22}S - S\tilde{A}_{22}^T + \tilde{B}_{22}R^{-1}\tilde{B}_{22}^T &= 0, \label{chap how far pert:math: Ss cond}
	\end{align}
	where $\tilde{A}_{22}$, $\tilde{B}_{22}$ and $\tilde{C}_{22}$ are obtained from the Schur decomposition (\ref{math: Schur decomp equations}). Let $Z$ be defined as in (\ref{math: schur Z definition}) and let $X = T-S$. Suppose, $\lambda \neq 0$ is an eigenvalue of $\tilde{A}_{22}$ with corresponding eigenvector ${v_2}^T \in ker(\tilde{A}_{22}^T - \lambda{I})$ and let ${y}=U[ 0 \quad  {v_2}^T]^T \in ker(A_r^T - \lambda{I})$, where $U$ is obtained from the Schur decomposition (\ref{math: Schur decomp equations}). Then ${y}^TZ{y} < 0$ if and only if ${v_2}^T X {v_2}>0$.
\end{lemma}
\begin{proof}
	Suppose all assumptions are satisfied and (\ref{chap how far pert:math: Ts cond})-(\ref{chap how far pert:math: Ss cond}) hold. Subtracting (\ref{chap how far pert:math: Ts cond}) from (\ref{chap how far pert:math: Ss cond}) gives the Lyapunov equation,
	\begin{align}
		0 &= -\tilde{A}_{22}S - S\tilde{A}_{22}^T + \tilde{B}_{22}R^{-1}\tilde{B}_{22}^T \nonumber \\
		& \hspace{80pt}  - (-\tilde{A}_{22}T - T\tilde{A}_{22}^T + \tilde{C}_{22}R\tilde{C}_{22}^T ) \nonumber\\
		&= \tilde{A}_{22}X + X\tilde{A}_{22}^T + \tilde{Z}_{22}. \label{chap how far pert:math: 22 Lyapunov equation}
	\end{align}
	Note that the matrix $\tilde{A}_{22}$ is the $22$ block of the matrix $\tilde{A}$ obtained from the Schur decomposition (\ref{math: Schur decomp equations}). By assumption, we know there exists a vector ${v_2}$ that satisfies
	\begin{align*}
		{v_2}^T\tilde{A}_{22} = \lambda{{v_2}^T}, \\
		\tilde{A}_{22}^T {v_2} = \bar{\lambda} {v_2}.
	\end{align*}
	Now, pre-multiply (\ref{chap how far pert:math: 22 Lyapunov equation}) by ${v_2}^T$ and post-multiply by ${v_2}$ to form
	\begin{align}
		{v_2}^T \tilde{A}_{22}X {v_2} + {v_2}^T X\tilde{A}_{22}^T {v_2} + {v_2}^T\tilde{Z}_{22}{v_2} = 0.
	\end{align}
	We may proceed with the following manipulation:
	\begin{align}
		0 &= {v_2}^T \tilde{A}_{22}X {v_2} + {v_2}^T X\tilde{A}_{22}^T {v_2} + {v_2}^T\tilde{Z}_{22}{v_2} \nonumber\\
		&= \lambda{{v_2}^T}X {v_2} + {v_2}^T X\bar{\lambda} {v_2} + {v_2}^T\tilde{Z}_{22}{v_2} \nonumber\\
		&= (\lambda + \bar{\lambda}){v_2}^T X {v_2} + {v_2}^T\tilde{Z}_{22}{v_2}. \label{math: lambda equation for x = z}
	\end{align}	
	We now rearrange (\ref{math: lambda equation for x = z}) to form
	\begin{align}	
		{v_2}^T X {v_2} &= - \frac{{v_2}^T\tilde{Z}_{22}{v_2}}{(\lambda + \bar{\lambda})}, \label{chap how far pert:math: X &= Zf_22 equation}
	\end{align}
	where, $(\lambda + \bar{\lambda}) > 0$. 
	In order to see the relationship between ${v_2}^T X {v_2}$ and ${y}^TZ{y}$, observe that
	\begin{align}
		{v}^T\tilde{Z}{v} &=\begin{bmatrix}
			0 & {v_2}^T
		\end{bmatrix}
		\begin{bmatrix}
			\tilde{Z}_{11} & \tilde{Z}_{12} \\ \tilde{Z}_{21} & \tilde{Z}_{22}
		\end{bmatrix}
		\begin{bmatrix}
			0 \\ {v_2}
		\end{bmatrix} = {v}_{2}^T\tilde{Z}_{22}{v_2}. \label{chap how far pert:math: vZfv = vZf22v}
	\end{align}
	Using (\ref{chap how far pert:math: vZfv = vZf22v}), we can rewrite (\ref{chap how far pert:math: X &= Zf_22 equation}) as
	\begin{align}
		{v_2}^T X {v_2} = - \frac{{v}^T\tilde{Z}{v}}{(\lambda + \bar{\lambda})}. \label{chap how far pert:sni: math: v2Xv2 = -Zf/lambda}
	\end{align}
	We know that $ {y}^T = {v}^TU^T$ by assumption. Thus we can substitute back into (\ref{chap how far pert:sni: math: v2Xv2 = -Zf/lambda}) to see that
	\begin{align}
		{v_2}^T X {v_2} &= - \frac{{v}^T\tilde{Z}{v}}{(\lambda + \bar{\lambda})} = - \frac{{y}^TZ{y}}{(\lambda + \bar{\lambda})}. \label{X, tilde Z, Z relationship eq}
	\end{align}
	Since $(\lambda + \bar{\lambda}) > 0$ and (\ref{X, tilde Z, Z relationship eq}) was arrived at using only algebraic manipulation, we can infer that ${y}^TZ{y} < 0$ if and only if ${v_2}^T X {v_2}>0$.
\end{proof}

\vspace{10pt}

\begin{lemma} \label{chap how far pert:Lemma: Z cases for X pos cor} 
	Consider a SISO  system (\ref{uncertain system}) satisfying Assumptions A1-A2. Let $Z,A_r$ and $X=T-S$ be defined as in (\ref{math: schur Z definition}), (\ref{A_r equation}) and (\ref{math: T-S cond}). Suppose there exists a vector ${y}$ such that ${y}=U[ 0 \quad  {v_2}]^T \in ker(A_r^T - \lambda{I})$, where $\lambda > 0$ and $U$ is obtained from the Schur decomposition (\ref{math: Schur decomp equations}). Under these assumptions, ${v_2}^TX{v_2}=0$ if and only if either of the following conditions holds
	\begin{align} 
		\frac{(C_1B_1B_2^T-C_1B_2B_1^T)}{(C_1B_2)^2}{y} &= 0, \label{ Z cases for X pos cor: cond: a}\\
		\text{or  } \qquad  B_2^T{y} &= 0. \label{ Z cases for X pos cor: cond: b}
	\end{align} 
\end{lemma}

\begin{proof}
	Suppose either (\ref{ Z cases for X pos cor: cond: a}) or (\ref{ Z cases for X pos cor: cond: b}) hold. We know from Lemma~\ref{chap how far pert:Lemma: Z cases for X pos def} that this immediately implies ${y}^T Z {y}=0$. Also, we know from Equation (\ref{X, tilde Z, Z relationship eq}) in Lemma~\ref{ch how far perturb: lem Z<0 implies X > 0 : non scalar} that for a vector ${y}$ such that ${y}=U[ 0 \quad  {v_2}]^T \in ker(A_r^T - \lambda{I})$, then
	\begin{align}
		{v_2}^T X {v_2} &= - \frac{{v}^T(\tilde{Z}){v}}{(\lambda + \bar{\lambda})} = - \frac{{y}^TZ{y}}{(\lambda + \bar{\lambda})}, \label{math: X = Z = 0}
	\end{align}
	where, $(\lambda + \bar{\lambda}) > 0$. Since $(\lambda + \bar{\lambda}) > 0$, it follows from (\ref{math: X = Z = 0}) that ${y}^T Z {y}=0$ implies ${v_2}^T X {v_2}=0$.
	
	Conversely, suppose ${v_2}^T X {v_2}=0$. We can use the relationship (\ref{math: X = Z = 0}) and after substituting in ${v_2}^T X {v_2}=0$ we see that
	\begin{align}
		- \frac{{y}^TZ{y}}{(\lambda + \bar{\lambda})} = 0.
	\end{align}
	Since, $(\lambda + \bar{\lambda}) > 0$ by assumption, this implies
	\begin{align}
		{y}^TZ{y} = 0.
	\end{align}
	Now, Lemma~\ref{chap how far pert:Lemma: Z cases for X pos def} states that ${y}^T Z {y}=0$ if and only if (\ref{ Z cases for X pos cor: cond: a})-(\ref{ Z cases for X pos cor: cond: b}) hold.
\end{proof}

\vspace{10pt}

\begin{lemma} \label{chap how far pert}
	Consider the  system (\ref{uncertain system}) satisfying Assumptions A1-A2. Let $Z,A_r$ and $X=T-S$ be defined as in (\ref{math: schur Z definition}), (\ref{A_r equation}) and (\ref{math: T-S cond}). Suppose the Schur decomposition (\ref{math: Schur decomp equations}) is such that $dim(\tilde{A}_{22})=1$, and let ${y}=U[ 0 \quad  1]^T \in ker(A_r^T - \lambda{I})$, where $\lambda = \tilde{A}_{22} > 0$ and $U$ is obtained from the Schur decomposition (\ref{math: Schur decomp equations}).
	Then $X=0$, if and only if either of the following conditions holds
	\begin{align} 
		\frac{(C_1B_1B_2^T-C_1B_2B_1^T)}{(C_1B_2)^2}{y} &= 0, \label{ Z cases for X pos cor: cond: a scalar}\\
		\text{or  } \qquad  B_2^T{y} &= 0. \label{ Z cases for X pos cor: cond: b scalar}
	\end{align} 
\end{lemma}

\begin{proof}
	This is a special case of Lemma~\ref{chap how far pert:Lemma: Z cases for X pos cor} and the proof is the same noting that ${v_2}$ in that proof is equal to $1$ in this case.
	We know from Lemma~\ref{ch how far perturb: lem Z<0 implies X > 0 : non scalar} that for a vector ${y}$ such that ${y}=U[ 0 \quad  {v_2}]^T \in ker(A_r^T - \lambda{I})$, $X$ is related to $Z$ by
	\begin{align}
		{v_2}^T X {v_2} &=  - \frac{{y}^TZ{y}}{(\lambda + \bar{\lambda})}. \label{math: scalar X proof: X Z relationship}
	\end{align}
	When $dim(\tilde{A}_{22})=1$, ${v_2} = 1$ and we can substitute ${y}$ into (\ref{math: scalar X proof: X Z relationship}) to obtain
	\begin{align}
		X  &= - \frac{\begin{bmatrix}
				0 & 1 \end{bmatrix}^T\tilde{Z}\begin{bmatrix}
				0 & 1 \end{bmatrix}}{(\lambda + \bar{\lambda})} = - \frac{{y}^TZ{y}}{(\lambda + \bar{\lambda})}, \label{math: scalar X = Z = 0}
	\end{align}
	where, $(\lambda + \bar{\lambda}) > 0$. We know from Lemma~\ref{chap how far pert:Lemma: Z cases for X pos def} that ${y}^T Z {y}=0$ if and only if (\ref{ Z cases for X pos cor: cond: a scalar})-(\ref{ Z cases for X pos cor: cond: b scalar}) hold. Since $(\lambda + \bar{\lambda}) > 0$, it follows from (\ref{math: scalar X = Z = 0}) that ${y}^T Z {y}=0$ if and only if $X=0$.
\end{proof}

\vspace{10pt}

\begin{lemma} \cite{HJ12} \label{lemma: different eigenvalue, vectors orthog}
	Let matrix $A$, scalars $\lambda$ and $\mu$, and non zero vectors ${x}$, ${y}$ be given. Suppose that $A{x}=\lambda {x}$ and ${y}^TA=\mu{y}^T$.  If $\lambda \neq \mu$, then ${y}^T{x} = 0$.
\end{lemma}

\vspace{10pt}

\begin{lemma} \cite{HJ12} \label{lemma: lambda multi 1 implies y x neq 0}
	Let matrix $A$, eigenvalue $\lambda$ and non zero vectors ${x}$ and ${y}$ be given. Suppose that $\lambda$ is an eigenvalue of $A$ such that $A{x}=\lambda {x}$ and ${y}^TAx=\lambda{y}^T$. 
	If $\lambda$ has algebraic multiplicity $1$, then ${y}^T{x} \neq 0$.
	Also, if $\lambda$ has geometric multiplicity $1$, then it has algebraic multiplicity $1$ if and only if ${y}^T{x} \neq 0$.
\end{lemma}

\vspace{10pt}

The following lemma is important as it tells us how the eigenvalues of the matrix $A_q$ behave under the perturbation $A_r = A_q + {\epsilon}Q$. Remember the eigenvalues of $A_q$ correspond to the zeros of the transfer function from $u$ to $z$ in (\ref{uncertain system}), plus a zero at the origin. Also, the solvability of the Lyapunov equations in Theorem~\ref{theorem: SNI ARE synthesis theorem} is also depend on the eigenvalues of $A_r$.

\vspace{10pt}

\begin{lemma} \label{chap how far pert:chap SNI state: Theorem : Perturbation properties of Ar vectors}
	Consider the uncertain system (\ref{uncertain system}) with $m=1$ satisfying Assumptions A1-A2. Consider the $n \times n$ matrix $A_q = QA$ defined in (\ref{math: A_q equation}). Suppose $A_q$ is not a defective matrix and has the distinct eigenvalues
	\begin{align}
		0 \, , \,  \lambda_1 \, , \, \lambda_2 \, , \, ... \, , \, \lambda_k,
	\end{align}
	such that
	\begin{align}
		dim \big( \, ker( A_q) \, \big) &= j , \label{chap how far pert:math: dim  ker( A_q)= j }\\
		dim \big( \, ker( A_q - \lambda_i {I} ) \, \big) &= n_i, i = 1,2,...,k, \label{chap how far pert:math: dim  ker( A_q - lambda_i )  = n_i}
	\end{align}
	where $1 \leq j \leq n$. Let a value of $\epsilon > 0$ be given such that $\lambda_i \neq - \epsilon \quad \forall \, i$ where $i = 1,2,...,k$ and the matrix $A_r = QA + {\epsilon}Q$ is not defective. Then $A_r = QA+\epsilon{Q}$ has the distinct possible eigenvalues
	\begin{align}
		0 \, , \, \epsilon \, , \, \lambda_1 + \epsilon \, , \, ... \, , \, \lambda_k + \epsilon,
	\end{align}
	such that
	\begin{align}
		dim \big( \, ker( A_r - \epsilon{I}) \, \big) &= j - 1  \label{chap how far pert:math: dim ker ar - ep cond}\\
		dim \big( \, ker( A_r - (\lambda_i + \epsilon){I}) \, \big) &= n_i, i = 1,2,...,k \\
		dim \big( \, ker( A_r ) \, \big) &= 1 \label{chap how far pert:math: dim ker ar = 1 cond}.
	\end{align}
\end{lemma}
\begin{proof} The proof will be published elsewhere. \end{proof}

\vspace{10pt}

\subsection{Main results}

In this section we address the problem of determining the maximum degree of stability that can be achieved using the SNI state feedback control laws presented in Theorem~\ref{theorem: SNI ARE synthesis theorem no A22} and Theorem~\ref{theorem: SNI ARE synthesis theorem}. These results are given under a number of assumptions on the eigenvalues of the matrix $A_q$. Future research will be directed towards more generalized results. Our results are offered as two separate cases that consider when the system has different zeros of the transfer function from $u$ to $z$ in (\ref{uncertain system}).

\vspace{10pt}

For the results that follow, we assume the following assumption is satisfied:

\vspace{10pt}

\begin{assumption}
$A_q$ has distinct eigenvalues $\lambda_1, \cdots , \lambda_{n}$. \label{ass 3}
\end{assumption}
This assumption is made so that we can leverage Lemma~\ref{chap how far pert:chap SNI state: Theorem : Perturbation properties of Ar vectors}, which does not use generalized eigenvector. This assumption will be relaxed in future work.

\vspace{10pt}

{\em Case 1:} 

The following case considers when the matrix $A_q$ has a single unstable eigenvalue, a single eigenvalue at the origin, and $n-2$ stable eigenvalues. We only consider a single unstable eigenvalue here in order to simply the proof and conserve space.

\vspace{10pt}

The eigenvalues of $A_q$ are ordered such that
\begin{align}
	Re[\lambda_1] \leq \cdots \leq Re[\lambda_{n-2}] \leq 0 \leq \lambda_n
\end{align}
and  we will assume that the following conditions hold:
\begin{enumerate}[(i)]
	\item For all $i \in \{1,\cdots,n-2\}$, $Re[\lambda_i] < 0$;
	\item For $i = n-1$, $\lambda_i= 0$;
	\item For $i = n$, $\lambda_i> 0$.
\end{enumerate}
Let $\gamma$ be defined as $\gamma = -Re[\lambda_{n-2}]$. In other words, $\gamma$ is the real part of the maximum eigenvalue of $A_q$ less than zero.

\vspace{10pt}

In the following theorem we will show that a system with an $A_q$ matrix as above may be perturbed for all $\epsilon$ within the range $0 < \epsilon < \gamma$ and still satisfy the conditions of Theorem \ref{theorem: SNI ARE synthesis theorem}.

\vspace{10pt}

\begin{theorem} \label{chap how far pert: distinct eigenvalues}
	Consider the system (\ref{uncertain system}) that satisfies the conditions of Theorem \ref{theorem: SNI ARE synthesis theorem}. Suppose, the pair $\{A,B_2\}$ is controllable,  the matrix $A_q$, (\ref{math: A_q equation}), has  distinct eigenvalues and the matrix $A_q$, has only one unstable eigenvalue;
	The closed-loop system formed using Theorem \ref{theorem: SNI ARE synthesis theorem} will be SNI with degree of stability $\epsilon$, for any $\epsilon$ within the range $0 < \epsilon < \gamma$, where $\gamma = -Re[\lambda_{n-2}]$.
\end{theorem}

\vspace{10pt}

{\em Case 2:} 

The following case considers when the matrix $A_q$ has no unstable eigenvalues, a single eigenvalue at the origin, and $n-1$ stable eigenvalues.

\vspace{10pt}

Now, consider a matrix $A_q \in \mathbb{R}^{n \times n}$ with distinct eigenvalues $\lambda_1, \cdots , \lambda_{n}$. Suppose the eigenvalues of $A_q$ are ordered such that
\begin{align}
	Re[\lambda_1] \leq \cdots \leq Re[\lambda_{n-1}] < 0. \label{math: gamma for no a22 case n x n}
\end{align}
Also, assume that the following conditions hold:
\begin{enumerate}[(i)]
	\item For all $i \in \{1,\cdots,n-1\}$, $Re[\lambda_i] < 0$;
	\item For $i = n$, $\lambda_i= 0$.
\end{enumerate}
 Let $\gamma$ be defined as $\gamma = -Re[\lambda_{n-1}]$. In other words, $\gamma$ is the real part of the maximum eigenvalue of $A_q$ less than zero.

\vspace{10pt}

In the following theorem we show that the system with an $A_q$ matrix as above may be perturbed for all $\epsilon$ within the range $0 < \epsilon < \gamma$ and still satisfy the conditions of Theorem \ref{theorem: SNI ARE synthesis theorem no A22}; that is, the system will be SNI for all $\epsilon$ within the range $0 < \epsilon < \gamma$. Moreover, for $\epsilon > \gamma$, the system may still be SNI if it satisfies Theorem~\ref{chap how far pert: distinct eigenvalues}.

\vspace{10pt}

\begin{theorem}  \label{chap how far pert: distinct eigenvalues no A22}
	Consider the system (\ref{uncertain system}) with dimension $n$, that satisfies the conditions of Theorem \ref{theorem: SNI ARE synthesis theorem no A22} when unperturbed. Suppose, 
	 the pair $\{A,B_2\}$ is controllable;
	 the matrix $A_q$ defined in  (\ref{math: A_q equation}) has distinct eigenvalues, and 
	the matrix $A_q$, has no unstable eigenvalues.
Then, the closed-loop system formed using Theorem \ref{theorem: SNI ARE synthesis theorem no A22} will be SNI with degree of stability $\epsilon$, for any $\epsilon$ within the range $0 < \epsilon < \gamma$, where $\gamma = -Re[\lambda_{n-1}]$ is defined as in (\ref{math: gamma for no a22 case n x n}). In addition,  for $\epsilon$ in the range $-Re[\lambda_{n-1}] < \epsilon < -Re[\lambda_{n-2}]$, if the open-loop system satisfies the conditions of Theorem \ref{theorem: SNI ARE synthesis theorem} at $\epsilon$, then the system will remain SNI for all $\epsilon$ within the range $-Re[\lambda_{n-2}] < \epsilon < -Re[\lambda_{n-2}]$.
\end{theorem}  

\vspace{10pt}

{\em Proof of Theorem \ref{chap how far pert: distinct eigenvalues}:}
	If a system satisfies the conditions of Theorem \ref{theorem: SNI ARE synthesis theorem}, there exists an $X>0$ such that
	\begin{align}
		\tilde{A}_{22}X + X\tilde{A}_{22}^T + \tilde{B}_{22}R^{-1}\tilde{B}_{22}^T - \tilde{C}_{22}R\tilde{C}_{22}^T = 0 \nonumber\\
		\tilde{A}_{22}X + X\tilde{A}_{22}^T + \tilde{Z}_{22} = 0, \label{ch pertb: n x n: gen layap}
	\end{align}
	where $\tilde{A}_{22}$, $\tilde{B}_{22}$ and $\tilde{C}_{22}$ are obtained from the Schur decomposition (\ref{math: Schur decomp equations}). Suppose the conditions of Theorem \ref{theorem: SNI ARE synthesis theorem} are satisfied when $\epsilon=0$ and $X_0>0$ is a solution to (\ref{ch pertb: n x n: gen layap}). The dimension of $X_0$ is related to $\tilde{A}_{22}$ by $dim(X_0) = dim(\tilde{A}_{22})=1$, so in this case $X_0$ is scalar. Our proof methodology will show that for any $\epsilon$ within the range $0 < \epsilon < -Re[\lambda_{n-2}]$, if there exists a left eigenvector of $A_r$ such that $X_\epsilon=0$, then $X_0=0$ must also be true, which contradicts our initial assumption. To that end, let $\gamma = -Re[\lambda_{n-2}]$ and suppose there exists an $\epsilon$ within the range $0 < \epsilon < \gamma$ such that $X=0$. Also, let ${y_1}$ be a left eigenvector of $A_r$ such that 
	\begin{align}
		{y_1}^TA_r = (\mu + \epsilon) {y_1}^T, \label{ch pertb: n x n thm: mu eq}
	\end{align}
	where $A_r$ is defined in (\ref{A_r equation}) and $\mu=\lambda_n$ is the anti-stable eigenvalue of $A_q$.
	
	According to Lemma~\ref{chap how far pert}, if $X=0$, then either
	\begin{align} \label{chap how far pert:conditions: nxn case:  SNI two conditions for X > 0}
		\frac{(C_1B_1B_2^T-C_1B_2B_1^T)}{(C_1B_2)^2}{y_1} &= 0, \\
		\text{or  } \qquad  B_2^T{y_1} &= 0. \label{chap how far pert:conditions: nxn : SNI two conditions for X > 0, cond b}
	\end{align}
	
	We now consider each of these conditions separately:
	
	{\em Case 1, (\ref{chap how far pert:conditions: nxn case:  SNI two conditions for X > 0}) is satisfied:} 
	\\\\
	In this case,
	\begin{align}
		(\frac{C_1B_1B_2^T-C_1B_2B_1^T}{(C_1B_2)^2}){y_1} = {w}^T{y_1} &= 0. \label{chap how far pert:math: nxn: w^tC_1^T = 0 equation}
	\end{align}
	Since we can easily check that the vector ${w}$ is orthogonal to $C_1$, it follows that
	\begin{align}
		C_1^T \in ker(A_r^T). 
	\end{align}	
	Moreover, we have distinct eigenvalues, therefore $C_1$ is also a basis for $ker(A_r^T)$. Since ${w}$ is orthogonal to the basis of $ker(A_r^T)$, this must mean ${w} \in range(A_r)$. Thus, we can describe ${w}$ using a linear combination of basis vectors for $range(A_r)$, to give
	\begin{align}
		{w} = \beta_1 {x_1} + \beta_2 {x_2} + \cdots + \beta_{n-1}{x_{n-1}}, \label{nxn lem: math: w definition}
	\end{align}
	where $\beta_1,\cdots,\beta_{n-1} \in \mathbb{R}$ and ${x_1},\cdots,{x_{n-1}}$ are eigenvectors of $A_r$ corresponding to the $n-1$ non-zero eigenvalues.
	
	Using this form of ${w}$, we can restate (\ref{chap how far pert:math: nxn: w^tC_1^T = 0 equation}) as
	\begin{align}
		{w}^T{y_1} &= 0 \\
		\iff {y_1}^T{w} &= 0 \\
		\iff \beta_1 {y_1}^T{x_1} + \cdots + \beta_{n-1} {y_1}^T{x_{n-1}} &= 0. \label{nxn lem: math: w condition}
	\end{align}
	According to Lemma~\ref{lemma: different eigenvalue, vectors orthog}, every term after $\beta_1 {y_1}^T{x_1}$ will be equal to zero. Therefore (\ref{nxn lem: math: w condition}) is reduced to
	\begin{align}
		\beta_1 {y_1}^T{x_1} = 0. \label{nxn lem: math: w condition reduced}
	\end{align}
	Equation (\ref{nxn lem: math: w condition reduced}) tells us that Condition (\ref{chap how far pert:math: nxn: w^tC_1^T = 0 equation}) holds if either $\beta_1=0$ or ${y_1}^T{x_1} = 0$.
	
	Since we only have simple eigenvalues, it follows from Lemma~\ref{lemma: lambda multi 1 implies y x neq 0} that ${y_1}^T{x_1} = 0$ cannot occur. Therefore, we now consider $\beta_1=0$. 
	
	If we substitute $\beta_1=0$ back into (\ref{nxn lem: math: w definition}), the first term is eliminated and we are left with
	\begin{align}
		{w} = \beta_2 {x_2} + \cdots + \beta_{n-1}{x_{n-1}}. \label{nxn lem: math: w mu2x2 definition}
	\end{align}
	We know from Lemma~\ref{chap how far pert:chap SNI state: Theorem : Perturbation properties of Ar vectors} that $A_q$ and $A_r$ will have the same right-eigenvectors corresponding to a particular eigenvalue. Therefore, each vector ${x_i}$ is an eigenvector of $A_q$ and the following holds:
	\begin{align}
		A_q {x_i} = \lambda_i {x_i},
	\end{align}
	for $i = 1,2,\cdots, n-2,n$.
	Also, for the eigenvalue $\mu=\lambda_n$ in (\ref{ch pertb: n x n thm: mu eq}), there will be a left-eigenvector of $A_q$, ${y_0}$, that satisfies
	\begin{align}
		{y_0}^T A_q = \mu {y_0}^T. 
	\end{align}
	Now, consider 
	\begin{align}
		{w}^T{y_0} = \beta_2 {x_2}^T{y_0} + \cdots + \beta_{n-1}{x_{n-1}}{y_0}.
	\end{align}
	We can apply Lemma~\ref{lemma: different eigenvalue, vectors orthog} to show that each term in this sum is orthogonal and therefore zero. Thus,
	\begin{align}
		{w}^T{y_0} = \beta_2 {x_2}^T{y_0} + \cdots + \beta_{n-1}{x_{n-1}}{y_0} = 0. \label{math: nxn: wv1 is orthog}
	\end{align}
	Equation (\ref{math: nxn: wv1 is orthog}) is just (\ref{ Z cases for X pos cor: cond: a}). Therefore, we can conclude from Lemma~\ref{chap how far pert:Lemma: Z cases for X pos cor} that (\ref{math: nxn: wv1 is orthog}) implies $X_0=0$. This contradicts our initial assumption that $X_0>0$.
	\\\\
	{\em Case 2, (\ref{chap how far pert:conditions: nxn : SNI two conditions for X > 0, cond b}) is satisfied:} 
	\\\\
	In this case, let ${y_1}$ be a left eigenvector of $A_r$ such that 
	\begin{align}
		{y_1}^TA_r = (\mu + \epsilon) {y_1}^T
	\end{align} and
	\begin{align}
		B_2^T{y_1} &= 0.  \label{chap how far pert:2nd nxn condition}
	\end{align}
	We know $(A,B_2)$ is controllable by assumption. Since is $(A,B_2)$ controllable, it follows immediately that $(A + \epsilon{I} , B_2)$ is also controllable. Using this knowledge, we will make a manipulation on $A_r$ to show that it is also controllable. To that end, remember $A_r = QA + \epsilon{Q}$, $Q= I - B_2(C_1B_2)^{-1}C_2$, and $C_1B_2$ is a scalar. Let $k = -\frac{C_1(A+\epsilon{I})}{C_1B_2}$ and consider the following manipulation:
	\begin{align}
		A_r = Q(A + \epsilon{I})
		&= (A+\epsilon{I}) - (\frac{B_2C_1}{C_1B_2})(A+\epsilon{I}) \nonumber\\
		&= A + \epsilon{I} + B_2K. \label{chap how far pert:math: nxn ar is controllable}
	\end{align}
	Thus, we can conclude from (\ref{chap how far pert:math: nxn ar is controllable}) that $(A_r, B_2)$ is also controllable. An immediate consequence of $(A_r, B_2)$ being controllable is  that there can be no left eigenvector of $A_r$ that is orthogonal to $B_2$. Therefore, (\ref{chap how far pert:2nd nxn condition}) can never be satisfied and we have a contradiction. 
	\\\\
	We have shown that both cases lead to a contradiction. 
\hfill $\blacksquare$

\vspace{10pt}

{\em Proof of Theorem \ref{chap how far pert: distinct eigenvalues no A22}:} 
	For $\epsilon$ within the range $0 < \epsilon < \gamma$, where $\gamma = -Re[\lambda_{n-1}]$, the closed-loop system is SNI follows directly from Theorem~\ref{theorem: SNI ARE synthesis theorem no A22}. For $-Re[\lambda_{n-2}] < \epsilon < -Re[\lambda_{n-2}]$, the matrix $A_q$ will have attained an unstable eigenvalue. Consequently, $X$ has dimension one. Therefore, if the system satisfies the conditions of Theorem \ref{theorem: SNI ARE synthesis theorem} for $-Re[\lambda_{n-2}] < \epsilon < -Re[\lambda_{n-2}]$, then the system will remain SNI for all $\epsilon$ within the range $0< \epsilon < -Re[\lambda_{n-2}]$ follows directly from Theorem \ref{chap how far pert: distinct eigenvalues}. 
\hfill $\blacksquare$

\vspace{10pt}

\section{Illustrative example} \label{chap how far pert:chap state feedback: worked example}

In this section we provide an illustrative example that illustrates Theorem~\ref{chap how far pert: distinct eigenvalues}.
 
\vspace{10pt}

Consider an uncertain system

\begin{align*}
	A &= \begin{bmatrix}
		-1 & 0 & -1 \\ 1 & 0 & -1 \\ -1 & 2 & 1
	\end{bmatrix}, &
	B_1 &= \begin{bmatrix}
		1 \\ 1 \\ 1
	\end{bmatrix}, &
	B_2 &= \begin{bmatrix}
		0 \\ 1 \\ 1
	\end{bmatrix}, \\
	C_1 &= \begin{bmatrix}
		1 & 1 & 0
	\end{bmatrix}, & D &= 0.
\end{align*}
that satisfies assumptions A1-A2. The $A_q$ matrix for this system is 
\begin{align*}
	A_q &= \begin{bmatrix}
		-1 & 0 & -1 \\ 
		1 & 0 & 1 \\ 
		-1 & 2 & 3
	\end{bmatrix},
\end{align*}
with $\sigma(A_q) = \{-1.6458 , \quad 0 , \quad 3.6458 \}$. We can readily verify that for any $\epsilon<1.6458$ this system satisfies the conditions of Theorem~\ref{theorem: SNI ARE synthesis theorem}, and the closed-loop system is SNI.

Now, we choose $\epsilon=1.6458$ and reevaluate. With $\epsilon=1.6458$, the matrix $A_r = QA+\epsilon{Q}$ has the form
\begin{align*}
	A_r &= \begin{bmatrix}
		$0.6458$ & $0$ & $-1$ \\ 
		$-0.6458$ & $0$ & $1$ \\ 
		$-2.6458$ & $0.3542$ & $4.6458$
	\end{bmatrix},
\end{align*}
with $\sigma(A_r) = \{0 , \quad 5.2946 \}$ and the $0$ eigenvalue has multiplicity two. At exactly $\epsilon=1.6458$, the matrix $X$ is undergoing a dimension change and cannot be evaluated. Therefore, we analyse $X$ as we approach $\epsilon=1.6458$ from the positive direction on the real axis; i.e., we consider $\epsilon^+ : \epsilon \to 1.6458^+$.

We can reform $X_{\epsilon^+}$ for our perturbed system as 	\begin{align*}
	X_{\epsilon^+} &= \begin{bmatrix}
		$-5.5431$ & $0.2815$ \\ 
		$0.2815$ & $0.2477$
	\end{bmatrix}.
\end{align*}
It is easily verified that $X_{\epsilon^+}$ is dimension 2 and no longer positive definite, with a negative eigenvalue at $-5.5568$. We no longer satisfy the conditions of Theorem \ref{theorem: SNI ARE synthesis theorem} and cannot guaranteed the SNI property.

\section{Conclusions}
In this paper we have studied the maximum prescribed degree of stability that can be achieved by SISO closed-loop systems formed using SNI state feedback. These results show that in the cases considered, the achievable degree of stability for the closed-loop system is related to the zeros of the transfer of the nominal plant from the control input to the disturbance output. These results were given under a number of assumptions and future research will be directed towards relaxing these assumptions and extending the results to MIMO systems. 

\addtolength{\textheight}{-12cm}   





\bibliography{irpnew}
\bibliographystyle{IEEEtran}

\end{document}